\documentclass[11pt]{amsart}

\usepackage{latexsym,amsthm,amssymb}

\usepackage[leqno]{amsmath}

\usepackage{mathrsfs}

\usepackage[all]{xy}

\newtheorem{theorem}{Theorem}[section]
\newtheorem{lemma}[theorem]{Lemma}
\newtheorem{definition}[theorem]{Definition}
\newtheorem{proposition}[theorem]{Proposition}

\newtheorem{remark}[theorem]{Remark}

\newtheorem{corollary}[theorem]{Corollary}

\numberwithin{equation}{theorem}

\newenvironment{Proof}
{\noindent \emph{Proof.}}
{\hfill $\square$ \medskip}

\newenvironment{Proofof}
{\noindent \emph{Proof of}}
{\hfill $\square$ \medskip}

\renewcommand{\mathcal}{\mathscr}

\newcommand{\Cal}{\mathcal}

\newcommand{\SE}{{\mathcal{E}}}
\newcommand{\SF}{{\mathcal{F}}}

\newcommand{\SH}{{\mathcal{H}}}
\newcommand{\SI}{{\mathcal{I}}}

\newcommand{\SL}{{\mathcal{L}}}

\newcommand{\SN}{{\mathcal{N}}}
\newcommand{\SO}{{\mathcal{O}}}

\newcommand{\SX}{{\mathcal{X}}}
\newcommand{\SY}{{\mathcal{Y}}}

\renewcommand{\mathbb}{\mathbf}

\newcommand{\PP}{\mathbb{P}}

\newcommand{\CC}{\mathbb{C}}

\newcommand{\WY}{\widetilde{Y}}
\newcommand{\WX}{\widetilde{X}}

\newcommand{\WL}{\widetilde{L}}

\newcommand{\wphi}{\widetilde{\varphi}}
\newcommand{\wrho}{\widetilde{\rho}}

\title{K3 double structures on Enriques surfaces and their smoothings}

\author{Francisco Javier Gallego}
\author{Miguel Gonz\'alez}
\author{\\ Bangere P. Purnaprajna}

\address{Departamento de \'Algebra, Universidad Complutense de Madrid}
\email{gallego@mat.ucm.es}
\address{Departamento de \'Algebra, Universidad Complutense de Madrid}
\email{mgonza@mat.ucm.es}
\address{Department of Mathematics, University of Kansas}
\email{purna@math.ku.edu}


\begin{document}
\begin{abstract}
Let $Y$ be a smooth Enriques surface. A $K3$ carpet on $Y$ is a locally Cohen-Macaulay
double structure on $Y$ with the same invariants as a smooth $K3$ surface (i.e., regular and with trivial canonical sheaf). The surface
$Y$ possesses an \'etale $K3$ double cover
$X \overset{\pi} \longrightarrow Y$.  We prove that $\pi$ can be deformed to a family $\SX \longrightarrow \mathbf P^N_{T^*}$ of projective embeddings of $K3$ surfaces and that any projective $K3$ carpet on $Y$ arises from such a family as the
flat limit of smooth, embedded $K3$ surfaces.
\end{abstract}
\maketitle
\section*{Introduction}

In this article we study
the relation between double covers and the origin of double structures. This relation was first studied  
in~\cite{Fong}, for hyperelliptic canonical morphisms and the so-called canonical ribbons and in~\cite{GP} for hyperelliptic $K3$ surfaces and $K3$ carpets on rational normal scrolls. Recently, M. Gonz\'alez in~\cite{Gon} and the authors in~\cite{GGP} studied this relation in a much more general setting, namely, finite covers of curves of arbitrary degree on the one hand and one dimensional, locally Cohen-Macaulay multiple structures of arbitrary multiplicity on the other hand.
In the present work we look at the relation between a natural and particularly nice double cover, the \'etale $K3$ double cover of an Enriques surface, and an interesting class of double structures, the {\it $K3$ carpets on Enriques surfaces}.

\medskip

Double structures on surfaces have appeared in connection with the study of the zero locus of sections of the Horrocks--Mumford vector bundle (see for example the work of Hulek and Van de Ven in~\cite{HV}). Also, $K3$ carpets on rational normal scrolls have been considered in the study of degenerations of smooth $K3$ surfaces. In this article we study another kind of $K3$ carpets, namely, those supported on Enriques surfaces.
A $K3$ carpet on a smooth Enriques surface $Y$ will be a locally Cohen-Macaulay
double structure on $Y$ with the same invariants as a smooth $K3$ surface (i.e., regular and with trivial canonical sheaf). The surface $Y$ possesses an \'etale $K3$ double cover
$X \overset{\pi} \longrightarrow Y$ associated to the canonical bundle of $Y$, which is $2$-torsion.
We prove that any projective $K3$ carpet on $Y$ arises from a family $\SX \longrightarrow \mathbf P^N_{T^*}$ of projective embeddings of $K3$ surfaces that degenerates to $\pi$.
As a consequence of this, we show that any projective $K3$ carpet on $Y$ can be {\it smoothed}, i.e., obtained as the flat limit of a family of smooth, irreducible (projective $K3$) surfaces.

\medskip

The reader might probably have noted in the previous paragraph the phrase ``projective $K3$ carpet".
$K3$ carpets on an Enriques surfaces (like indeed double structures on any other surface) need not be projective, unlike ribbons on curves. Thus our first task is to characterize (see Theorem~\ref{nonproj.carpets}) those $K3$ carpets which are projective. This is accomplished in Section~\ref{embedding}.  There we also see ``how many" projective $K3$ carpets there are. We do this in two settings. On the one hand, we compare the sizes of the families of
projective $K3$ surfaces on a given (abstract) Enriques surfaces $Y$ and the size of the family of non--projective $K3$ carpets (see Theorem~\ref{nonproj.carpets}). This situation has some strong resemblance to
the case of projective and non--projective smooth $K3$ surfaces, where the former lie on infinite, countably many codimension $1$ families in the moduli space of $K3$ surfaces. On the other hand we also compute the dimension of the space that parametrizes the family of projective $K3$ surfaces supported on a given Enriques surface which is embedded in a projective space (see Theorem~\ref{dim.K3carpet.embedded}).

\medskip

In Section~\ref{smoothing} we prove the results regarding deformation of morphisms and smoothings of carpets. First we show (see Theorem~\ref{coversmoothing}) that the cover $\pi$ can be deformed to a family of embeddings of $K3$ surfaces to projective space. Then, in order to obtain a smoothing of a projective $K3$ carpet $\WY$, one considers a suitable embedding of $\WY$ in projective space, then one chooses the family of embeddings of Theorem~\ref{coversmoothing} suitably, in order to obtain a family of projective schemes consisting of the images of smooth $K3$ surfaces degenerating to $\WY$. From these theorems we obtain a smoothing result for most of the embedded $K3$ carpets (see Theorem~\ref{embsmoothing}) and subsequently we show that any (abstract) projective $K3$ carpet can be smoothed (see Theorem~\ref{absgensmooth}).

\medskip

Finally we devote Section~\ref{hilbert.section} to study the Hilbert points of projective $K3$ carpets. We prove that their Hilbert point are always smooth (see Theorem~\ref{hilbert}), unlike the case of $K3$ carpets on rational normal scrolls (in that case, some Hilbert points are smooth and some are not; see~\cite[Section 4]{GP}).

\medskip
\noindent
{\bf Acknowledgements}: We thank Joseph Lipman for a helpful discussion and for pointing out some references regarding the dualizing sheaf. We also thank N. Mohan Kumar for some useful discussions.

\medskip
\noindent
{\bf Convention.} {\rm We work over $\CC$. 
Throughout the article, when we talk about a regular or an Enriques surface, we will mean it to be smooth, irreducible and proper over $\CC$.}

\section{$K3$ carpets. Characterization}\label{characterization}
Among carpets on an Enriques surface $Y$, we single out a family which deserve special attention as far as they share the invariants of smooth $K3$ surfaces. We call them $K3$ carpets.
In fact, we will give a more general definition: a $K3$ carpet on any regular surface will be a carpet with the same invariants of a smooth $K3$ surface (i.e., trivial dualizing sheaf and irregularity $q=0$; see Definition~\ref{K3carpet} and Proposition~\ref{char2.K3}).
Gallego and Purnaprajna, in~\cite{GP}, studied $K3$ carpets supported on rational normal scrolls.
In this paper we consider carpets on a different type of surfaces possessing a double covering from a smooth $K3$ surface, namely Enriques surfaces.
In this new case, as in~\cite{GP}, the adjective $K3$ is not only justified by the fact that these carpets have the same invariants as smooth $K3$ surfaces, but also from the fact that
projective $K3$ carpets are degenerations of smooth $K3$ surfaces, as we shall prove in this paper.

\medskip

\noindent
We start by recalling the definition of a carpet on a smooth surface.

\noindent
\begin{definition}\label{defrope}
Let $Y$ be a reduced connected scheme and let $\SE$ be a line bundle on $Y$.
A ribbon on $Y$ with conormal bundle $\SE$ is a scheme $\WY$ with ${\WY}_{\mathrm{red}}=Y,$ such that
\begin{enumerate}
\item
$\SI_{Y, \WY}^2=0$ and
\item
$\SI_{Y, \WY} \simeq \SE$ as $\SO_{Y}$--modules.
\end{enumerate}
When $Y$ is a surface, $\WY$ is called a carpet on $Y$.
\end{definition}

\noindent
We give now the definition of a $K3$ carpet supported on a regular surface.  Although our definition does not require the carpet to be a regular scheme, we will see in Proposition~\ref{char2.K3} that a $K3$ carpet defined according to Definition~\ref{K3carpet} is always regular.
\begin{definition}\label{K3carpet}
Let $Y$ be a regular surface. A $K3$ carpet $\WY$ on $Y$ is a carpet on $Y$ such that its dualizing sheaf $\, \omega_{\WY}\simeq \SO_{\WY}$.
\end{definition}

\noindent
The existence of a dualizing sheaf with nice functorial properties on a proper scheme is not obvious.
In Remark~\ref{dualizing} we justify the existence of the dualizing sheaf in Definition~\ref{K3carpet}.
In Lemma~\ref{dualizing.properties} we point out some nice properties of the dualizing sheaf on $\WY$.
The assertions in Remark~\ref{dualizing} and in Lemma~\ref{dualizing.properties} are valid, in general, for ribbons.
\begin{remark}\label{dualizing}
{\rm Let $Y$ be a smooth irreducible proper variety.
\begin{enumerate}
\item Any ribbon $\WY$ on $Y$ is a proper scheme over $\CC$.
So, according to, e.g.,~\cite[(7), p. 46]{Kleiman80}, there is a dualizing sheaf $\, \omega_{\WY}$ on $\WY$.
\item Any ribbon $\WY$ on $Y$ is a locally Gorenstein (in fact, locally a complete intersection) scheme. Therefore the dualizing sheaf $\, \omega_{\WY}$ is an invertible sheaf, see \cite[V 9.3, 9.7, VII 3.4]{RD} and \cite[p. 157]{Conrad00}.{\hfill $\square$ \medskip}
\end{enumerate}
}
\end{remark}

\begin{lemma}\label{dualizing.properties}
Let $Y$ be a smooth irreducible proper variety.
Let $\WY$ be a ribbon on $Y$ with conormal bundle $\SE$ and $\omega_{\WY}$ its dualizing sheaf.
\begin{enumerate}
\item The dualizing sheaf on $Y$ is
\begin{equation*}\label{dualizing1}
 {\SH}om_{\WY}(\SO_Y, \omega_{\WY})= \omega_Y.
\end{equation*}

\item Let $\SL$ be an invertible sheaf on $Y$.
Then
\begin{equation*}\label{dualizing2}
{\SH}om_{\WY}(\SL, \omega_{\WY})= \SL^{-1}\otimes \omega_Y.
\end{equation*}

\item The dualizing sheaf on $\WY$ fits into an extension
\begin{equation}\label{WY-extension.dual.2}
\xymatrix@1@C-5pt{
0 \ar[r] & \omega_{Y} \ar[r] & \omega_{\WY} \ar[r] & \SE^{-1} \otimes \omega_Y  \ar[r] & 0.}
\end{equation}
\end{enumerate}
\end{lemma}
\begin{proof}
(1) From the definition of a dualizing sheaf, see~\cite[p. 241]{RD} or~\cite[(1),(6)]{Kleiman80},
we see that ${\SH}om_{\WY}(\SO_Y, \omega_{\WY})$ is a dualizing sheaf on $Y$.\\\\
\noindent
(2) Since $\SL$ is a sheaf on $Y$, notice that ${\SH}om_{\WY}(\SL, \omega_{\WY})={\SH}om_{Y}(\SL, {\SH}om_{\WY}(\SO_{Y}, \omega_{\WY}))$.\\\\
\noindent
(3) From the inclusion $\SE \hookrightarrow \SO_{\WY}$, we have a map
$\omega_{\WY} \to {\SH}om_{\WY}(\SE, \omega_{\WY})$.
We see at once that, since $\omega_{\WY}$ is invertible, this map is surjective.
So applying $\mathcal{H}om_{\WY}(-,\omega_{\WY})$ to
\begin{equation}\label{WY-extension.2}
\xymatrix@1@C-5pt{
0 \ar[r] & \SE \ar[r] & \SO_{\WY} \ar[r] & \SO_{Y} \ar[r] & 0,}
\end{equation}
we obtain~\eqref{WY-extension.dual.2}.
\end{proof}

\medskip
\noindent
Now we characterize $K3$ carpets from its conormal bundle.
\begin{proposition}\label{char.K3}
Let $Y$ be a regular surface and let $\WY$ be a carpet whose reduced part is $Y$.
Let $\SE$ be the ideal sheaf of $Y$ in $\WY$. Then $\WY$ is a $K3$ carpet iff $\,\SE \simeq \omega_Y$.
\end{proposition}
\begin{Proof}
Let $\mathcal{E} \simeq \omega_{Y}$.
Look at~\eqref{WY-extension.dual.2}.
Since $H^1(\omega_Y)=0$,  the section $1 \in H^0(\SO_Y)$ can be lifted to $H^0(\omega_{\WY})$, and hence, $\omega_{\WY}$ being invertible, we have $\omega_{\WY} \simeq \SO_{\WY}$.
Now assume $\, \omega_{\WY}\simeq \SO_{\WY}$.
If we tensor~\eqref{WY-extension.dual.2} with $\SO_Y$
we get a surjection $\SO_Y \to \SE^{-1} \otimes \omega_Y$. Thus $\SE \simeq \omega_Y$.
\end{Proof}

\noindent
As a consequence of~\ref{char.K3} we see that a $K3$ carpet, as defined in Definition~\ref{K3carpet}, is a regular scheme, as is the case of smooth $K3$ surfaces.
\begin{proposition}\label{char2.K3}
Let $\WY$ be a $K3$ carpet on a regular surface $Y$. Then $H^{1}(\SO_{\WY})=0$.
\end{proposition}
\begin{Proof}
From Proposition~\ref{char.K3}, the conormal bundle $\SE = \omega_Y$.
Since $Y$ is a regular surface $H^{1}(\omega_{Y})=H^{1}(\SO_{Y})=0$ and hence, from~\eqref{WY-extension.2}, $H^{1}(\SO_{\WY})=0$.
\end{Proof}

\begin{remark}\label{Ext1}
{\rm We have seen that the $K3$ carpets on a given  regular surface are the carpets with conormal bundle $\omega_Y$.
Thus (see~\cite[1.4]{BE95}) the space of non--split $K3$ carpets on a given regular surface $Y$ is the projective space of lines in $\mathrm{Ext}_{Y}^1(\Omega_Y, \omega_Y)$.

\noindent
Notice that, when $Y$ is an Enriques surface the dimension of $\mathrm{Ext}_{Y}^1(\Omega_Y, \omega_Y)$ is the Hodge number $h^{1,1}=10$.{\hfill $\square$ \medskip}}
\end{remark}

\section{Projective and non--projective $K3$ carpets}\label{embedding}

\noindent
In contrast to ribbons on curves, not all carpets are
projective, (see ~\cite[III Exercise
5.9]{Hartshorne77}) even if all of them are proper
or if, as is the case with Enriques surfaces, they
are supported on a projective surface.
Thus the very first question about the $K3$ carpets on Enriques
surfaces is whether there exist families of projective $K3$ carpets. This question
has a positive answer as is illustrated in Theorem 2.2 and Theorem 2.3.  Next step
is to compute the dimension of the space parametrizing $K3$ carpets on a given Enriques surface.
This is settled in Theorem~\ref{dim.K3carpet.embedded} for the dimension of the family of
embedded (projective) carpets on a given embedded Enriques surface,
and in Theorem ~\ref{nonproj.carpets}, where we compute the size of
the space of projective $K3$ carpets supported on a given (abstract) Enriques
surface $Y$, comparing it also with the space of all $K3$ carpets on $Y$.
As we will see, the situation somehow resembles
that of smooth $K3$ surfaces.

\smallskip

\noindent
To start searching for embedded $K3$ carpets we need to look first for
embeddings of Enriques surfaces in projective space. We recall some
well known facts about this:

\begin{remark}\label{Enriques.embeddings}
{\rm Let $Y$ be an Enriques surface.
\begin{enumerate}
\item
If $Y$ is embedded in $\mathbf P^N$, then $N \geq 5$.
\item
A very ample line bundle on $Y$ has sectional genus $g \geq 6$ and degree $d \geq 10$.
\item
If $N \geq 5$, then the surface $Y$ can be embedded in $\mathbf P^N$.
\end{enumerate}}
\end{remark}

\begin{proof}
By adjunction, there do not exist Enriques
surfaces in $\mathbf P^3$.
On the other hand, applying the formula for the numerical
invariants of a smooth surface $Y$ in $\mathbf P^4$
(see~\cite[A.4.1.3]{Hartshorne77}),
\begin{equation*}
d^{2}-10d-5H K_{Y}-2K_{Y}^2+12+12p_{a}=0,
\end{equation*}
we see at once that there do not exist Enriques surfaces in $\mathbf
P^4$ either. This completes the proof of (1).
Now, a line bundle on $Y$ with sectional genus $g$ has $g$ linearly
independent global section. Then, if the line bundle is very ample,
(1) implies that $g \geq 6$, so its
degree is $2g-2 \geq 10$. This proves (2). Finally,
since $Y$ is projective, $Y$ can be embedded in $\mathbf P^M$,
with $M >> 0$ and we project it isomorphically into $\mathbf P^N$ as far as $N \geq 5$.
\end{proof}

\noindent
Now we want to know how many $K3$ carpets are supported on a given embedded
Enriques surfaces. This will do in Theorem~\ref{dim.K3carpet.embedded}. To do this we will need to know the
dimension of the space of first--order infinitesimal deformations of a morphism from a $K3$ surface to projective space.
Given a morphism $\varphi$ from a variety $X$ to $\PP^N$, the normal sheaf $\SN_{\varphi}$ is defined as the cokernel of the natural map $\Cal{T}_{X}\longrightarrow  \varphi^{\ast} \Cal{T}_{\PP^N}$. Then the first--order infinitesimal deformations of $\varphi$, up to isomorphism, are parametrized by $H^0(\SN_{\varphi})$ (see~\cite[4.2]{Hor}).
In our setting since $X$ is a smooth $K3$ surface, it is a smooth variety. Then, if the image of $\varphi$ has the same dimension as $X$, we have the following exact sequence:

\begin{equation}\label{tan.normal.varphi.1}
0 \longrightarrow \Cal{T}_{X}\longrightarrow  \varphi^{\ast} \Cal{T}_{\PP^N}
\longrightarrow \SN_{\varphi} \longrightarrow 0.
\end{equation}

\begin{theorem}\label{dim.19}
Let $X$ be a smooth projective $K3$ surface and let $X \xrightarrow{\varphi} \PP^N$
be a morphism whose image is a surface.
Let $\SN_{\varphi}$ be the normal sheaf of $\varphi$. Then,
\begin{enumerate}
\item
the dimension of the image of the connecting map
\begin{equation*}
H^0(\SN_{\varphi}) \to H^1(\Cal{T}_X)
\end{equation*}
of the long exact sequence of cohomology of~\eqref{tan.normal.varphi.1}
is $19$;
\item
$H^1(\SN_{\varphi})=0$; and
\item $H^2(\SN_{\varphi})=0$.
\end{enumerate}
\end{theorem}

\begin{Proof}
Let us denote $L=\varphi^{\ast}\SO_{\PP^N}(1)$ and let us consider the Atiyah extension of $L$
\begin{equation}\label{Atiyah.L}
0 \longrightarrow \SO_X \longrightarrow \Sigma_{L} \longrightarrow \Cal{T}_{X}\longrightarrow 0.
\end{equation}
The space $H^1(\Sigma_L)$ parametrizes first--order infinitesimal deformations of the pair $(X, L)$ up to isomorphism (see \cite[pp. 126--128]{Zar} or \cite[II.2.2]{Sernesi06}) and
the map $H^0(\SN_{\varphi}) \to H^1(\Cal{T}_X)$ factors through $H^1(\Sigma_L)$.
Taking cohomology on~\eqref{Atiyah.L} yields the exact sequence
\begin{equation*}
H^1(\Cal O_X) \longrightarrow  H^1(\Sigma_L) \longrightarrow  H^1(\Cal{T}_{X}) \longrightarrow H^2(\Cal O_X).
\end{equation*}
Since $X$ is a $K3$ surface, $h^1(\Cal O_X)=0$, $h^2(\Cal O_X)=1$ and
$h^1(\Cal{T}_{X})$ is the same as the Hodge number
$h^{1,1}$ of $X$, hence
\begin{equation}\label{dim4}
\dim H^{1}(\Cal{T}_{X})=20.
\end{equation}
On the other hand,  $H^{1}(\Cal{T}_X) \longrightarrow  H^2(\SO_X) $ is induced by cup product with the cohomology class $c(L) \in H^1(\Omega_X)$ (see \cite[Proposition II.2.2]{Sernesi06}),  so it is surjective, for $L$ is non--trivial (see \cite[p. 57]{Sernesi06}).
Then
\begin{equation}\label{dim.SigmaL}
\dim H^{1}(\Sigma_L) = 19.
\end{equation}
Then, going back to ~\eqref{tan.normal.varphi.1} we have the long exact sequence
\begin{equation*}
H^{0}(\SN_{\varphi}) \xrightarrow{\nu}
H^{1}(\Cal{T}_X) \longrightarrow
H^{1}(\varphi^{\ast}\Cal{T}_{\PP^{N}}) \longrightarrow H^{1}(\SN_{\varphi}) \longrightarrow 0,
\end{equation*}
where the exactness on the far right comes from $h^2(\Cal{T}_X)=h^{0,1}=0$.
Then~\eqref{dim.SigmaL} implies that the image of $\nu$ has dimension less than or equal to $19$. On the other hand, taking cohomology on the dual of the Euler sequence restricted to $X$ yields $h^1(\varphi^{\ast}\Cal{T}_{\PP^{N}})=1$, for $H^1(L)=H^2(L)=0$ since $L$ is ample. All this together with~\eqref{dim4} implies that the image of $\nu$ has dimension  $19$ and $H^{1}(\SN_{\varphi})=0$.

\noindent
To prove (3) note that taking cohomology on the dual of the Euler sequence restricted to $X$ yields $H^2(\varphi^{\ast}\Cal{T}_{\PP^{N}})=0$, for
$H^2(L)=0$. Then it follows that $H^{2}(\SN_{\varphi})=0$.
\end{Proof}

\noindent
We will use Theorem~\ref{dim.19} in this situation (see e.g. ~\eqref{h1.normal.zero} in the proof of Theorem~\ref{dim.K3carpet.embedded}):
we set $\varphi$ to be the composition of the  \'etale $K3$ double cover $X \xrightarrow{\pi} Y$ of an Enriques surface $Y$ followed by an embedding $Y \overset{i}\hookrightarrow \PP^N$.
On the other hand, Theorem~\ref{dim.19} can be also used if $\varphi$ is an embedding into projective space, so we recover the following result:

\begin{corollary}\label{smooth.Hilbert}
If $X$ is a smooth projective $K3$ surface embedded in projective space, (not necessarily as a linearly normal variety nor as a non--degenerate variety), then the point of $X$ in the Hilbert scheme is smooth.
\end{corollary}

\noindent
Next theorem gives a quantitative measure on the $K3$ carpets supported on an embedded Enriques surfaces.
Precisely, given an embedded Enriques surface $Y\overset {i}{\hookrightarrow}\mathbf{P}^{N}$ , we find the dimension
of the variety that parametrizes the $K3$ carpets in $\mathbf P^N$, supported on $i(Y)$.

\begin{theorem}\label{dim.K3carpet.embedded}
Let $Y$ be an Enriques surface and let $Y \overset{i}{\hookrightarrow}
\mathbf P^N$ be an embedding of $Y$.
Let $g$ be the sectional genus of $i(Y)$.
The $K3$ carpets embedded in $\mathbf{P}^{N}$ and supported on
$i(Y)$ are parametrized by a non--empty open set in the projective space of lines in
$H^0(\SN_{Y,\PP^{N}} \otimes \omega_Y)$, whose dimension is
$g(N+1)+8$.
In particular, if $i$ is induced by the complete linear series of
$\Cal O_Y(1)$, then the dimension of 
this open set is $g^2+8$.
\end{theorem}

\begin{Proof}
Denote $\mathcal I=\mathcal{I}_{i(Y),\PP^N}$.
The $K3$ carpets in $\mathbf{P}^{N}$ which are supported on $i(Y)$
are in one--to--one correspondence with the surjective elements in
$\mathrm{Hom} \,(\mathcal{I}/\mathcal{I}^{2},\mathcal{\omega}_{Y})$,
up to nonzero scalar multiple (see ~\cite[Proposition
2.1.(2)]{Gon}; see also ~\cite[Lemma 1.4]{GP} or~\cite{HV}).

\noindent
We start computing the dimension of $\mathrm{Hom}
\,(\mathcal{I}/\mathcal{I}^{2},\mathcal{\omega}_{Y})$.
Recall that $\Omega_Y^{\ast} \otimes \omega_Y \simeq \Omega_Y$. Then,
since $Y$ is regular, and by Serre duality and Hodge Theory,
we have $h^0(\Omega_Y^{\ast} \otimes \omega_Y)= h^2(\Omega_Y^{\ast}
\otimes \omega_Y)=0$. Then, taking cohomology on the conormal sequence
of $i(Y)$, we get
\begin{multline}\label{hom.conormal.seq}
0 \to   \mathrm{Hom}\,(\Omega_{\PP^N} \otimes \SO_Y,\omega_Y) \to
\mathrm{Hom}\,(\mathcal{I}/\mathcal{I}^{2},{\omega}_{Y})
\overset{\delta} \longrightarrow
\mathrm{Ext}^{1}(\Omega_{Y},{\omega}_{Y}) \to  \\
\to \mathrm{Ext}^{1}(\Omega_{\PP^N}\otimes \SO_Y, \omega_{Y}) \to
 \mathrm{Ext}^{1}(\mathcal{I}/\mathcal{I}^{2},{\omega}_{Y}) \to 0.
\end{multline}
To find the dimension of  $\mathrm{Hom}
\,(\mathcal{I}/\mathcal{I}^{2},\mathcal{\omega}_{Y})$ we need to
compute the dimensions of the other terms of the sequence~\eqref{hom.conormal.seq}.
Dualizing the restriction to $Y$ of the Euler sequence and tensoring by
$\omega_{Y}$, we have the exact sequence
\begin{equation}\label{euler2}
\xymatrix@1@C-5pt{
0 \ar[r] & \omega_{Y} \ar[r] &
\mathcal{O}_Y^{\oplus
  N+1}(1) \otimes \omega_{Y} \ar[r] & \Omega_{\PP^N}^{\ast}
\otimes\omega_{Y} \ar[r] & 0.}
\end{equation}
Since
$h^{1}(\mathcal{O}_{Y}(1)\otimes\omega_{Y})=h^{2}(\mathcal{O}_{Y}(1)\otimes\omega_{Y})=0$,
it follows that $h^{1}(\Omega_{\PP^N}^{\ast} \otimes
\omega_{Y})=h^{2}(\omega_{Y})=1$.
So
\begin{equation}\label{dim1}
\dim \, \mathrm{Ext}^{1}(\Omega_{\PP^{N}} \otimes \SO_Y, {\omega}_{Y})=
h^{1}(\Omega_{\PP^{N}}^{\ast}\otimes {\omega}_{Y})=1.
\end{equation}

\noindent
Also, $h^{0}(\omega_{Y})=h^{1}(\omega_{Y})=0$, so we have
\begin{equation}\label{dim2}
\dim \, \mathrm{Hom}\,(\Omega_{\PP^{N}}\otimes \SO_Y,\omega_{Y})=
(N+1) \cdot h^{0}(\SO_{Y}(1)\otimes\omega_{Y})=g(N+1).
\end{equation}

\noindent
On the other hand (see Remark ~\ref{Ext1})
\begin{equation}\label{dim3}
\dim \, \mathrm{Ext}^{1}(\Omega_{Y},\omega_{Y})=10.
\end{equation}

\noindent
Finally we will see that $\mathrm{Ext}^{1}(\mathcal{I}/\mathcal{I}^{2}, \omega_{Y})=0$.
To do this, let
$X \overset{\pi}\rightarrow Y$ be the \'etale $K3$ double cover of $Y$.
Denote $i \circ \pi=\varphi$.
From Theorem~\ref{dim.19}, (2), for the normal sheaf of $\varphi$ we have
\begin{equation}\label{h1.normal.zero}
H^{1}(\SN_{\varphi})=0.
\end{equation}

\noindent
We will see that $\mathrm{Ext}^{1}(\mathcal{I}/\mathcal{I}^{2}, \omega_{Y})$ is a direct summand of $H^{1}(\SN_{\varphi})$.

\noindent Let $\SF$ be the kernel of
$\varphi^{\ast}\Omega_{\PP^N}\rightarrow\Omega_{X}$.
Since $\pi$ is \'etale, it follows that $\Omega_{X/Y}$ and
$\Omega_{X/\PP^N}$ are both $0$,
so we have the following commutative diagram:
\begin{equation*}\label{diagram}
\xymatrix@C-10pt@R-13pt{
& 0 \ar[d] &   0 \ar[d] &  &  \\
 & \pi^{\ast}(\mathcal{I}/\mathcal{I}^2) \ar[d] \ar@{=}[r] &
 \pi^{\ast}(\mathcal{I}/\mathcal{I}^2) \ar[d] &  & \\
0 \ar[r] & \SF \ar[d] \ar[r] & \varphi^{\ast}\Omega_{\PP^N} \ar[d] \ar[r] & \Omega_X \ar@{=}[d] \ar[r]& 0 \\
 & 0 \ar[r] & \pi^{\ast}\Omega_Y \ar[d] \ar[r] & \Omega_X \ar[r]& 0 \\
 & &  0.& &}
\end{equation*}

\noindent
Therefore there is an isomorphism
\begin{equation*}
\SN_{\varphi}\simeq \mathcal{H}om
\,(\pi^{\ast}\mathcal{I}/\mathcal{I}^{2},\mathcal{O}_{X}).
\end{equation*}
Since $\pi_{\ast}\mathcal{O}_{X}=\mathcal{O}_{Y}\oplus \omega_{Y}$,
taking cohomology and using the adjunction isomorphism we get
\begin{multline*}\label{H1normalphi}
H^{1}(\SN_{\varphi})= H^1(\mathcal{H}om \,(\pi^{\ast}\mathcal{I}/\mathcal{I}^{2},\mathcal{O}_{X}))=
\mathrm{Ext}^{1}(\pi^{\ast} \mathcal{I}/\mathcal{I}^{2}, \SO_{X})= \\
= \mathrm{Ext}^{1}(\mathcal{I}/\mathcal{I}^{2}, \SO_{Y})
\oplus \mathrm{Ext}^{1}(\mathcal{I}/\mathcal{I}^{2}, \omega_{Y}).
\end{multline*}
Then Theorem~\ref{dim.19}, (2) implies
\begin{equation}\label{dim6}
\mathrm{Ext}^{1}(\mathcal{I}/\mathcal{I}^{2}, \omega_{Y})=0.
\end{equation}

\noindent
Then, from~\eqref{hom.conormal.seq},~\eqref{dim1},~\eqref{dim2},~\eqref{dim3}, and~\eqref{dim6},
we see at once that
\begin{equation*}\label{dimension}
\dim \mathrm{Hom}\,(\mathcal{I}/\mathcal{I}^{2}, \omega_{Y})=g(N+1)+9.
\end{equation*}

\noindent
Recall that the $K3$ carpets on $Y$ embedded in $\mathbf{P}^{N}$ are
in one--to--one correspondence with the surjective homomorphisms in
$\mathrm{Hom}\, (\mathcal{I}/\mathcal{I}^{2},\mathcal{\omega}_{Y})$,
up to nonzero scalar multiple, or equivalently, with the nowhere
vanishing global sections of the $(N-2)$--rank vector bundle $\SN_{Y,
  \PP^N} \otimes \omega_Y$, up to nonzero scalar multiple.
Recall also that the elements of $\mathrm{Hom}\,
(\mathcal{I}/\mathcal{I}^{2},\mathcal{\omega}_{Y})$ corresponding to
surjective homomorphisms form an open set (see ~\cite[Lemma4.1]{Gon}).
Therefore, to finish the proof we need to show that there is a nowhere vanishing
section in the space $H^0(\SN_{Y, \PP^N} \otimes \omega_Y)$.
Observe first that $\SN_{Y,\PP^{N}}\otimes \omega_{Y}$ is globally generated.
To see this note that we have a surjection
$\Omega_{\PP^{N}}^{\ast} \otimes \omega_{Y}
\rightarrow \SN_{Y,\PP^{N}} \otimes \omega_{Y}$
so, from~\eqref{euler2},
we see that $\SN_{Y,\PP^{N}} \otimes \omega_{Y}$ is globally generated
as long as $\SO_{Y}(1)\otimes\omega_{Y}$ is globally generated.
This follows from Reider's theorem (\cite{Reider}), since
$\SO_{Y}(1)$ is very ample and its degree $d = 2g-2 \geq 10$ (see
Remark ~\ref{Enriques.embeddings}, (2)).
Finally since the rank of $\SN_{Y,\PP^{N}}\otimes
\omega_{Y}$ is $N-2 > \dim \,Y$ (see Remark~\ref{Enriques.embeddings}, (1))
and it is a globally generated vector bundle,
it has a nowhere vanishing section.
Thus the $K3$ carpets  inside $\PP^{N}$, supported on $i(Y)$ in $\PP^{N}$,
are parametrized by a non--empty open set in the projective space of
lines in $H^0(\SN_{Y,\PP^{N}} \otimes \omega_Y)$, whose dimension is
$g(N+1)+8$.
\end{Proof}

\smallskip
\noindent
The following theorem is a refinement of  ~\cite[III Ex. 5.9]{Hartshorne77}
to characterize non--projective $K3$ carpets.
As result of this theorem, we can say more about the size of  the families of
projective $K3$ carpets on a given (abstract) Enriques surface,
compared to the set of non--projective $K3$ carpets.

\begin{theorem}\label{nonproj.carpets}
Let $Y$ be an Enriques surface and let $\WY$ be a
$K3$ carpet on $Y$ corresponding to an element $\tau \in
\mathrm{Ext}^1(\Omega_Y,\omega_Y)$.
\begin{enumerate}
\item The carpet $\WY$ is projective if
and only if there exists an ample divisor $D$ on $Y$ such that $\int_D
\tau =0$, when $\tau$ is thought as an element of $H^{1,1}(Y)=H^2(Y,
\mathbf C)$.
\item
Non--split projective $K3$ carpets on $Y$
are parametrized by a union of (countably infinitely many distinct)
hyperplanes of
the $9$-dimensional projective space of lines in
$\mathrm{Ext}^1(\Omega_Y,\omega_Y)$. These hyperplanes are
in one--to--one correspondence with the set of
classes in NS$(Y)$ of primitive ample divisors on $Y$.
\end{enumerate}
\end{theorem}

\begin{proof}
Recall (see Remark ~\ref{Ext1}) that a $K3$ carpet on $Y$
corresponds to an element
\begin{equation*}
\tau \in \mathrm{Ext}^1(\Omega_Y,\omega_Y) \simeq H^1(\Omega_Y^*
\otimes \omega_Y) \simeq H^1(\Omega_Y)=H^{1,1}(Y)=H^2(Y,
\mathbf C).
\end{equation*}
Since the ideal of $Y$ inside $\WY$ is a square zero ideal,
we have an exact sequence
\begin{equation*}\label{like.exp}
0 \rightarrow \omega_{Y} \rightarrow \SO_{\WY}^{\ast} \rightarrow
\SO_{Y}^{\ast} \rightarrow 1.
\end{equation*}
This yields
\begin{equation*}
0 \rightarrow \mathrm{Pic}\, \WY \overset{\gamma}\longrightarrow
\mathrm{Pic}\, Y \overset{\lambda}\longrightarrow
H^{2}(\omega_{Y}) \rightarrow H^{2}(\SO_{\WY}^{\ast}) \rightarrow
H^{2}(\SO_{Y}^{\ast}).
\end{equation*}
The map $\lambda$ works as follows: if $D$ is  a divisor on $Y$, then
$\lambda(\Cal O_Y(D))=\int_D \tau$.
The map
$\gamma$ sends each line bundle on $\WY$ to its restriction to
$Y$. The carpet $\WY$ is projective if and only if it possesses an
ample line bundle.
On the other hand, a line bundle on $\WY$ is ample if and only
if its restriction to $Y$ is ample. Therefore $\WY$ is projective if
and only if there exists an ample line bundle on $Y$ that can be
lifted by $\gamma$ to $\WY$.
This is the same as saying that there
exists an ample line bundle on $Y$ lying in the kernel of $\lambda$.
Thus $\WY$ is projective if and only if there exists an ample divisor $D$
on $Y$ such that $\int_D \tau=0$.
Then, given an ample divisor $D$ on $Y$, the elements $\tau \in
\mathrm{Ext}^1(\Omega_Y,\omega_Y) \simeq H^{1,1}(Y)$ with $\int_D
\tau=0$ form a hyperplane $H_D$ of $\mathrm{Ext}^1(\Omega_Y,\omega_Y)$,
whose elements correspond to projective $K3$ carpets. Then projective
$K3$ carpets are parametrized by the projective lines in
\begin{equation*}
\bigcup_D H_D,
\end{equation*}
where $D$ ranges over the set of primitive ample divisors on $Y$.
\end{proof}

\begin{remark} {\rm Let $Y$ be an Enriques surface.
Theorem~\ref{nonproj.carpets} shows in particular the existence of
non--projective $K3$ carpets on a given Enriques surface $Y$.
Indeed, the non--split non--projective $K3$ carpets on $Y$ are parametrized by the
complement of a union of countably many hyperplanes of the
$9$-dimensional projective space of lines in
$\mathrm{Ext}^1(\Omega_Y,\omega_Y)$.
There are ``more'' non--projective
$K3$ carpets than projective $K3$ carpets.}
\end{remark}

\noindent
The arguments of the proof of Theorem ~\ref{dim.K3carpet.embedded} give
another way of looking at Theorem~\ref{nonproj.carpets}:

\begin{proposition}
Let $Y$ be an Enriques surface. Associated to every
embedding $i$ of $Y$ into some projective space $\mathbf P^N$, there is
a sequence~\eqref{hom.conormal.seq}, arising from the conormal
sequence of $i(Y)$ in $\mathbf P^N$.
For the sequence~\eqref{hom.conormal.seq}
associated to $i$, we will denote by  $\delta_i$ the
map $\delta$.
Let $\mathbf{P}(\mathrm{Im}\,\delta_i)$ be the projective space of lines in
$\mathrm{Im}\, \delta_i$.
Then the non--split projective $K3$ carpets on $Y$ are parametrized by
\begin{equation*}
\bigcup_i \mathbf P(\mathrm{Im}\, \delta_i),
\end{equation*}
where $i$ ranges among all the embeddings of $Y$ into some projective space.
For each $i$, $\mathbf P(\mathrm{Im} \, \delta_i)$ is a hyperplane in the
$9$-dimensional projective space of lines in $\mathrm{Ext}^1(\Omega_Y,\omega_Y)$.
\end{proposition}

\begin{proof}
If a $K3$ carpet $\WY$ on $Y$ is projective, it can be embedded in
  some projective space $\mathbf P^N$
by the complete linear series of a very ample
  line bundle.
This embedding induces
  an embedding $i$ of $Y$ as (a degenerate) subvariety of $\mathbf
  P^N$. Let $\Cal I$
  be the ideal sheaf of $Y$ in $\mathbf P^N$. Then, the carpet $\WY$
  embedded in $\mathbf P^N$ corresponds to an element of
  $\mathrm{Hom}(\Cal I/\Cal I^2, \omega_Y)$. Thus, the $K3$ carpet
  $\WY$, considered as an abstract scheme, corresponds to a point
  lying in the image of the map $\delta_i$. From ~\eqref{dim1} and
  ~\eqref{dim6} we gather
that the cokernel
  of $\delta_i$ has dimension $1$, hence the image of $\delta_i$ in
  $\mathrm{Ext}^1(\Omega_Y,\omega_Y)$ is a hyperplane. Thus the class
  in  $\mathrm{Ext}^1(\Omega_Y,\omega_Y)$
  of every
  projective $K3$ carpet lies in the image of the map $\delta_i$
  associated to some embedding $i$ of $Y$ into some projective space. Since
  obviously the classes lying in the image of any of the maps $\delta_i$
  correspond to projective $K3$ carpets, we see that non--split projective $K3$
  are parametrized by
\begin{equation*}
\bigcup_i \mathbf{P}(\mathrm{Im}\, \delta_i)
\end{equation*}
where $i$ ranges among all the embeddings of $Y$ into some projective space.
\end{proof}

\medskip
\noindent
In Theorem ~\ref{dim.K3carpet.embedded} and Theorem~\ref{nonproj.carpets}
we saw how many projective
$K3$ carpets there are supported on an Enriques surface.
In the next observation, we describe how embeddings by a complete linear
series of a $K3$ carpet look like.

\begin{remark}\label{K3carpet.embeddings}
{\rm Let $Y$ be an Enriques surface and let $\WY$ be a
projective $K3$ carpet on $Y$.
Assume that $\WY$ is embedded, as a non--degenerate subscheme into some
projective space, by the complete linear series of a very ample line
bundle.
Let $g$ be the sectional genus of $\SO_Y(1)=\SO_{\WY}(1) \otimes
\SO_Y$.
Then, from $H^1( \SO_Y(1)\otimes\omega_Y )=0$ and the exact sequence
\begin{equation*}
\xymatrix@1{
0 \ar[r]  & \omega_Y(1)  \ar[r] & \SO_{\WY}(1) \ar[r] & \SO_Y(1) \ar[r] & 0,}
\end{equation*}
we have
\begin{equation*}
H^0(\SO_{\WY}(1))=H^0(\SO_Y(1))\oplus H^0( \omega_Y(1) ).
\end{equation*}
Therefore the embedding induced on $Y$ is also given by the complete
linear series of $\SO_Y(1)$ and there is a diagram
\begin{equation*}\label{complete.embedding}
\xymatrix@C+5pt{
{\WY \;} \ar@{^{(}->}[r] & {\;\PP^{2g-1}=\PP(H^0(\SO_Y(1))\oplus
  H^0(\omega_Y(1)))}\\
{Y \;} \ar@{^{(}->}[u] \ar@{^{(}->}[r] & {\;\PP^{g-1}=\PP(H^0(\SO_Y(1))).} \ar@{^{(}->}[u]}
\end{equation*}
}
\end{remark}

\section{Deformation of morphisms and smoothing of projective
$K3$ carpets}\label{smoothing}

\noindent In this section we prove two results.
First we show in Theorem ~\ref{coversmoothing} that the  \'etale
$K3$ double cover $\pi$ of an Enriques
surface can be deformed, in many different ways, to a family of
projective embeddings. Second, as a consequence of Theorem
~\ref{coversmoothing} we show (see Theorem ~\ref{embsmoothing} and
Corollary ~\ref{absgensmooth}) that
every projective
$K3$ carpet $\WY$ on an
Enriques surface can be smoothed.
By this we mean that we can find a flat, proper, integral family $\SY$
over a smooth affine curve $T$, such that over for $0 \in T$,
$\SY_0=\WY$ and for $t \in T, t \neq 0$, $\SY_t$ is a
smooth, irreducible, and, in our case, projective $K3$ surface.

\medskip
\noindent
The key point that connects Theorems~\ref{coversmoothing}
and~\ref{embsmoothing} is the fact that $\WY$,
after being embedded in
some  projective space $\PP^N$, arises as the central fiber
of the image of a first--order infinitesimal deformation
of the composition of $\pi$ with the inclusion of $Y$ in $\PP^N$:

\begin{theorem}\label{infsmoothing}
Let $\, \WY \subset \PP^N$ be a projective $K3$ carpet on a smooth
Enriques surface $Y$. Let $X \overset{\pi}\to Y$ be the \'etale $K3$ double cover of $Y$
and let $X\overset{\varphi}\to \PP^N$ be the morphism obtained by composing
$\pi$ with the inclusion of $\,Y$ in $\PP^N$.
Then $\WY$ is the central fiber of the image of some first--order
infinitesimal deformation of $\,\varphi$.
\end{theorem}
\begin{Proof}
Since $\pi$ is \'etale, we have $\SN_{\pi}=0$. Then the result follows
from~\cite[Theorem 3.9]{Gon}.
\end{Proof}

\noindent
Next we show that $\varphi$ can be deformed to a
family of embeddings to $\mathbf P^N$. We do so by proving something
stronger, namely, that any infinitesimal
deformation of ${\varphi}$ can be extended to a family of embeddings
of smooth $K3$ surfaces in $\PP^N$. Theorem~\ref{coversmoothing} is,
in the present setting, the counterpart of ~\cite[Theorem 2.1]{GGP},
where the authors showed that
a finite cover of a curve can be deformed to a family of embeddings.

\begin{theorem}\label{coversmoothing}
Let $X  \overset{\pi} \longrightarrow Y$ be the \'etale $K3$ double
cover of an Enriques surface $Y$, embedded in $\mathbf P^N$ with
sectional genus $g$ and satisfying $N \leq 2g-1$. Let
$\varphi$ denote the composition of $\pi$ with the inclusion of $Y$ in
$\mathbf P^N$.
Let $\Delta=\text{Spec}\ {\mathbf k[\epsilon]}/{\epsilon^2}$. Then for
every first--order infinitesimal deformation
\begin{equation*}
\WX \overset{\wphi}\longrightarrow \PP_{\Delta}^N
\end{equation*}
of $X \overset{\varphi}\longrightarrow \PP^N$, there exists
 a smooth irreducible family $\SX$, proper and flat over a smooth
 pointed affine curve $(T, 0)$, and a $T$--morphism $\SX
 \overset{\Phi}\longrightarrow \PP_T^{N}$ with the following features:
\begin{enumerate}
\item the general fiber $\SX_t \overset{\Phi_t}\longrightarrow \PP^N$,
  $\, t \in T-0,$ is a closed immersion of a smooth $K3$ surface; and
\item the fiber of $\SX \overset{\Phi}\longrightarrow \PP_T^{N}$ over
  the tangent vector at $0 \in T$ is $\WX
  \overset{\wphi}\longrightarrow \PP_{\Delta}^{N}$; in particular, the
  central fiber $\SX_0 \overset{\Phi_0}\longrightarrow \PP^{N}$ is $X
  \overset{\varphi}\longrightarrow \PP^{N}$
\end{enumerate}
\end{theorem}

\begin{remark}
{\rm We require $N \leq 2g-1$ in the statement of
Theorem~\ref{coversmoothing}. This hypothesis is, in fact, quite natural.
Indeed, if $\, \WY \subset \PP^N$ is non--degenerate (i.e., not contained in a hyperplane),
then $N \leq 2g-1$ (see Remark~\ref{K3carpet.embeddings}).
The hypothesis is used in Step 2 of the proof
of Theorem~\ref{coversmoothing} (see~\eqref{bound}).}
\end{remark}

\noindent
Before proving Theorem~\ref{coversmoothing} we need the following
lemma:

\begin{lemma}\label{veryample}\noindent
Let $Y$ be an Enriques surface, embedded in projective space
  with sectional genus $g$,
 and let $X \overset{\pi}
  \longrightarrow Y$ be its \'etale $K3$ double cover.
Then, if   $L=\pi^{\ast}\SO_Y(1)$, $L$ is very ample.
\end{lemma}

\begin{Proof}
From  Remark ~\ref{Enriques.embeddings} it follows that $L^2=4g-4 \geq 20$.
Then, to prove that $L$ is very ample, it suffices to check the
following (see \cite[4.2, 5.2, 6.1]{SaintDonat74}):
\begin{enumerate}
\item there is no irreducible curve $E$ such that $p_a(E)=1$ and $L
  \cdot E=2$, and
\item there is no smooth rational curve $E$ such that $L \cdot E=0$.
\end{enumerate}
The first condition holds because $L$ is base-point-free and the
second condition holds because $L$ is ample.
\end{Proof}

\begin{Proofof}\emph{Theorem ~\ref{coversmoothing}.}
{Step 1.} To obtain $\Phi$ we first construct, in a suitable way, a
pair
$(\SX,\SL)$, where $\SX$ is a family of smooth $K3$ surfaces and $\SL$
is a family of very ample line bundles.\\
Let us denote
$\WL=\wphi^{\ast}\SO_{\PP_{\Delta}^{N}}(1)$.
Then $\WL$ restricts to $L$ on $X$ and
the $\Delta$--module $\Gamma(\WL)$ is free of rank $h^0(L)$ and
$\Gamma(\WL)\otimes k[\epsilon]/\epsilon k[\epsilon]= H^0(L)$. \\
Now we want to obtain a family $(\SX,\SL)$, proper and flat over a
smooth pointed affine curve $(T, 0)$, whose central fiber is $(X,L)$,
whose restriction to the tangent vector to $T$ at $0$ is $(\WX,\WL)$
and whose general member $(\SX_t,\SL_t)$ consists of a smooth
irreducible $K3$ surface and a very ample line bundle $\SL_t$.\\
Note that $L$ has degree $4g-4$ and $h^0(L)=2g$.
Then, from Lemma~\ref{veryample} we know that $L$ is very ample and,
by Corollary~\ref{smooth.Hilbert}, its
complete linear series $|L|$  defines an embedding
which determines a smooth point $[X]$ in a single
component of the Hilbert scheme of surfaces of degree $4g-4$ in $\mathbf
P^{2g-1}$.
The general point $[X']$ in this component represents a smooth
irreducible $K3$ surface.
Then we may consider an open neighborhood $H$ of $[X]$ in its Hilbert
component, with $H$ parametrizing only smooth $K3$ surfaces.
Moreover, since $L$ is very ample and $H^1(L)=0$,
also $\WL$ is very ample relative to $\Delta$ and the embedding $X
\overset{|L|}\hookrightarrow \PP^{2g-1}$ extends to an embedding $\WX
\hookrightarrow \PP_{\Delta}^{2g-1}$.
So the image of $\WX \hookrightarrow \PP_{\Delta}^{2g-1}$ is a flat
family over $\Delta$ which corresponds to a tangent vector to $H$ at
$[X]$.
We can take the embedding $\WX \hookrightarrow \PP_{\Delta}^{2g-1}$ so
that this tangent vector is nonzero.
Now, since $[X]$ is a smooth point in $H$, we can take a smooth
irreducible affine curve $T$ in $H$ passing through $[X]$ with tangent
direction the given tangent vector.\\
Let $0 \in T$ denote the point corresponding to $[X]$.
Then the pullback to $T$ of the universal family provides a family
$(\SX,\SL)$, proper and flat over $T$, whose central fiber is $(X,
L)$, whose restriction to the tangent vector to $T$ at $0$ is
$(\WX,\WL)$ and whose general member $(\SX_t, \SL_t)$ consists of a
smooth irreducible $K3$ surface and a very ample line
bundle $\SL_t$,
with $H^1(\SL_t)=H^2(\SL_t)=0$,
and hence, with $h^0(\SL_t)=h^0(L)=2g$.

\medskip
\noindent{Step 2.} Once we have the pair $(\SX,\SL)$, we are going to use it to construct a relative morphism
\begin{equation*}
\SX
\overset{\Phi}\longrightarrow \PP_T^{N}
\end{equation*}
with the properties described in the statement.

\medskip
\noindent Recall that $\SL$ is very ample relative to $T$ and that $h^0(\SL_t)=h^0(L)=2g$ and $h^1(\SL_t)=0$ for all $t \in T$.
Then formation of $p_*$ commutes with base extension and, after shrinking $T$, we may assume that $\Gamma(\SL)$ is a free $\SO_T$--module.
Then $\SL$ induces a morphism
 \begin{equation*}
\SX
\overset{\Psi}\longrightarrow \PP_T^{2g-1}
\end{equation*}
which is a closed immersion at each fiber.
The morphism $\wphi$ is the composition $\wrho \circ \Psi_\Delta$, for some linear projection
$\PP_{\Delta}^{2g-1} \overset{\wrho}\dashrightarrow \PP_{\Delta}^N$.
Now we look at some $t$ near (but different from) $0$. Since
\begin{equation}\label{bound}
N \leq 2g-1 =\dim \, |\SL_t|,
\end{equation}
we can find a linear projection $\rho_t$ mapping  $\Psi_t(\SX_t)$ to $\mathbf P^N$.
On the other hand, Remark ~\ref{Enriques.embeddings} implies $N \geq
5$. Then choosing $\rho_t$ sufficiently general, we may assume the composition $\rho_t \circ \Psi_t$ to be a closed immersion.
We lift $\wrho$ and $\rho_t$ to a linear projection $\rho$ to $\PP^N_T$.
Finally we define $\Phi$ as the composition $\rho \circ \Psi$. Since the restriction $\Phi_t$ is a closed immersion, by~\cite[4.6.7]{EGA3-1} so are the restrictions of $\Phi$ to the nearby fibers. Then, maybe shrinking $T$ we can conclude that the restriction of $\Phi$ to $\Delta$ is $\wphi$ and that the restrictions $\Phi_t$ are closed immersions for all $t \in T$, $t \neq 0$. \end{Proofof}

\medskip

\noindent
Now we use Theorems~\ref{infsmoothing} and~\ref{coversmoothing} to
show that $\WY$ is the limit of the images of a family of embeddings
$\Phi_t$ of smooth $K3$ surfaces, degenerating to
$\varphi$. Precisely, we want to extend the infinitesimal deformation
of $\varphi$ in such a way that, if we call the image of the family of
morphisms $\SY \subset \PP^N \times T$, then $\SY_0=\WY$. All this is
done in the next theorem:

\begin{theorem}\label{embsmoothing}
Let $\WY$ be a projective $K3$ carpet embedded in $\PP^N$, and
supported on an Enriques surface $Y$ embedded in $\mathbf P^N$
with sectional genus $g$ and $N \leq 2g-1$.
Then there exists a family of morphisms $\Phi$ over an affine curve
$T$ as described in Theorem~\ref{coversmoothing} such that the image
$\SY$ of $\Phi$  is a closed integral subscheme $\SY \subset \PP_T^N$,
flat over $T$, with the following features:
\begin{enumerate}
\item the general fiber $\,\SY_t, \, t \in T-0,$ is a smooth
  irreducible projective non--degenerate $K3$ surface in $\PP^N$,
\item the central fiber $\,\SY_0 \subset \PP^N$ is $\,\WY \subset \PP^N$.
\end{enumerate}
\end{theorem}
\begin{Proof}
We use the notations of the proof of the Theorem~\ref{coversmoothing}.\\
From Theorem~\ref{infsmoothing} we know that there exists a first
order infinitesimal deformation
\begin{equation*}
\WX \overset{\wphi}\to \PP_{\Delta}^{N}
\end{equation*}
of $\varphi$ such that the central fiber of the image of $\wphi$ is
equal to $\WY$.\\
Therefore there is a family $\SX \to T$ and a $T$--morphism $\SX \overset{\Phi}\to \PP_T^{N}$ as in Theorem~\ref{coversmoothing}.\\
Let $\SY$ be the image of the $T$--morphism $\SX \overset{\Phi}\to \PP_T^{N}$.
The total family $\SX$ is smooth and irreducible so $\SY$ is integral.
Furthermore, $\Phi$ is a closed immersion over $T-0$ since,
by~Theorem~\ref{coversmoothing}, $\Phi_t$ is a closed immersion for
every $t \in T-0$ (see e.g. \cite[4.6.7]{EGA3-1}).
Therefore for $t \in T-0$ we have the equality $\SY_t =
\,\mathrm{im}\,(\Phi_t)$. Since $\SX_t$ is smooth, this proves (1).
Finally, the fact that $T$ is an integral smooth curve and $\SY$ is
integral and dominates $T$ implies that $\SY$ is flat over $T$.
So the fiber $\SY_0$ of $\SY$ at $0 \in T$ is the flat limit of the
images of $\SX_t \overset{\Phi_t} \to \PP^{N}$ for $t \neq 0$.
Moreover, this fiber $\SY_0$ contains the central fiber
$\,(\mathrm{im}\,\widetilde{\varphi})_0$ of the image of $\wphi$.
Since $\WY$ has conormal bundle $\SE$ and $\pi$ has trace zero module
$\SE$, both $\SY_0$ and $(\mathrm{im}\,\widetilde{\varphi})_0$ have
the same Hilbert polynomial, so they are equal.
\end{Proof}

\medskip

We highlight this consequence of Theorem~\ref{embsmoothing}:

\begin{theorem}\label{absgensmooth}
Any projective $K3$ carpet $\, \WY$ on a Enriques surface $Y$
is smoothable.
\end{theorem}

\begin{proof}
Let us embed $\WY$ in projective space by the complete linear series
of a very ample line bundle.
Then Remark~\ref{K3carpet.embeddings} implies
that the condition $N \leq 2g-1$ is satisfied, so
the result follows from
Theorem~\ref{embsmoothing}.
\end{proof}

\section{The Hilbert point of a projective $K3$ carpet}\label{hilbert.section}
In this section we prove, in Theorem~\ref{hilbert}, that the Hilbert point of a projective $K3$ carpet on an Enriques surface is smooth. This is in sharp contrast with the result on Hilbert points corresponding to
$K3$ carpets on a rational normal scroll, as shown in  [GP97].
First we state some preliminary results valid in general for ribbons.
\begin{lemma}\label{dualizante-WY}
Let $\WY$ be a ribbon on a smooth irreducible proper variety $Y$ with conormal bundle $\SE$.
There is an isomorphism
\begin{equation}\label{omega.WY.restricted}
{\omega_{\WY} |}_{Y}=\SE^{-1} \otimes  \omega_{Y}.
\end{equation}
\end{lemma}
\begin{Proof}
Restricting the sequence~\eqref{WY-extension.dual.2}
to $Y$ gives the isomorphism.
\end{Proof}
\begin{lemma}\label{Hilbert}
Let $\, Y \subset \WY \subset \PP^N$ be an embedded ribbon, with conormal bundle $\SE$, on a smooth irreducible projective variety $Y$.
Then there are exact sequences
\begin{equation}\label{first-sequence-bis-twisted}
\xymatrix@1{
0 \ar[r] & {\SN_{\WY,\PP^N} |}_{Y}\otimes \SE \ar[r] & \SN_{\WY,\PP^N} \ar[r] & {\SN_{\WY,\PP^N} |}_{Y} \ar[r] & 0,}
\end{equation}
\begin{equation*}\label{seq4.}
\xymatrix@1{
0 \ar[r] &  \SE^{-1} \ar[r] & \SN_{Y, \PP^N} \ar[r]&  \mathcal{H}om_{Y}(\SI_{\WY,\PP^N}/\SI_{Y,\PP^N}^2, \SO_Y) \ar[r] & 0,}
\end{equation*}
and
\begin{equation*}\label{seq1.}
\xymatrix@1{
0 \ar[r] & \mathcal{H}om_{Y}(\SI_{\WY,\PP^N}/\SI_{Y,\PP^N}^2, \SO_Y)  \ar[r] & {\SN_{\WY, \PP^N} |}_{Y} \ar[r]& \SE^{-2} \ar[r] & 0.}
\end{equation*}
\end{lemma}
\begin{Proof}
We know that $\WY$ is a local complete intersection so $\SN_{\WY,\PP^N}$ is locally free. Therefore from
\begin{equation*}\label{first-sequence-bis}
\xymatrix@1{
0 \ar[r] & \SE \ar[r] & \SO_{\WY} \ar[r] & \SO_Y \ar[r] & 0,}
\end{equation*}
we obtain the sequence~\eqref{first-sequence-bis-twisted}.
Also $\SI_{\WY, \PP^N}/\SI_{\WY,\PP^N}^2$ is locally free so we have
\begin{equation*}
{\mathcal{H}om_{\WY}(\SI_{\WY, \PP^N}/\SI_{\WY,\PP^N}^2,\SO_{\WY}) |}_{Y}= \mathcal{H}om_{Y}({\SI_{\WY, \PP^N}/\SI_{\WY,\PP^N}^2 |}_{Y}, \SO_Y).
\end{equation*}
Furthermore ${\SI_{\WY, \PP^N}/\SI_{\WY,\PP^N}^2 |}_{Y}= \SI_{\WY,\PP^N}/\SI_{Y,\PP^N}\SI_{\WY,\PP^N}$ so we have an exact sequence
\begin{equation*}
\xymatrix@1{
0 \ar[r] & (\SE')^{-1} \ar[r] & {\SI_{\WY, \PP^N}/\SI_{\WY,\PP^N}^2 |}_{Y} \ar[r]& \SI_{\WY,\PP^N}/\SI_{Y,\PP^N}^2 \ar[r] & 0,}
\end{equation*}
where $\SE'$ is an invertible sheaf on $Y$.
So there is an exact sequence
\begin{equation}\label{seq1}
\xymatrix@1{
0 \ar[r] & \mathcal{H}om_{Y}(\SI_{\WY,\PP^N}/\SI_{Y,\PP^N}^2, \SO_Y)  \ar[r] & {\SN_{\WY, \PP^N} |}_{Y} \ar[r]& \SE' \ar[r] & 0,}
\end{equation}
Furthermore from
\begin{equation*}\label{seq3}
\xymatrix@1{
0 \ar[r] & \SI_{\WY,\PP^N}/\SI_{Y,\PP^N}^2  \ar[r] & \SI_{Y,\PP^N}/\SI_{Y,\PP^N}^2 \ar[r]& \SE \ar[r] & 0,}
\end{equation*}
we obtain the exact sequence
\begin{equation}\label{seq4}
\xymatrix@1{
0 \ar[r] &  \SE^{-1} \ar[r] & \SN_{Y, \PP^N} \ar[r]&  \mathcal{H}om_{Y}(\SI_{\WY,\PP^N}/\SI_{Y,\PP^N}^2, \SO_Y) \ar[r] & 0,}
\end{equation}
Moreover, since $\WY$ is a local complete intersection, we have
\begin{equation*}\label{seq6}
\bigwedge^{c}\SN_{\WY, \PP^N} =\omega_{\WY}\otimes \omega_{\PP^N}^{-1}= \omega_{\WY}\otimes \SO_{\WY}(N+1),
\end{equation*}
where $c$ is the codimension of $Y$.\\
So
\begin{equation*}\label{seq7}
{\bigwedge^{c}\SN_{\WY, \PP^N} |}_Y ={ \omega_{\WY} |}_Y \otimes \SO_Y(N+1),
\end{equation*}
and from the isomorphism~\eqref{omega.WY.restricted}
\begin{equation*}\label{seq8}
{\bigwedge^{c}\SN_{\WY, \PP^N} |}_Y =\omega_Y \otimes \SE^{-1} \otimes \SO_Y(N+1).
\end{equation*}
Moreover
\begin{equation*}\label{seq9}
\bigwedge^{c}\SN_{Y, \PP^N} =\omega_{Y}\otimes \SO_{Y}(N+1),
\end{equation*}
so
\begin{equation}\label{seq10}
{\bigwedge^{c}\SN_{\WY, \PP^N} |}_Y =\bigwedge^{c}\SN_{Y, \PP^N} \otimes \SE^{-1}.
\end{equation}
I claim that
\begin{equation}\label{SE'=SE-2}
\SE'=\SE^{-2}.
\end{equation}
Indeed, from~\eqref{seq1} we obtain
\begin{equation*}\label{seq11}
{\bigwedge^{c}\SN_{\WY, \PP^N} |}_Y =\bigwedge^{c-1} \mathcal{H}om_{Y}(\SI_{\WY,\PP^N}/\SI_{Y,\PP^N}^2, \SO_Y)\otimes \SE',
\end{equation*}
and from~\eqref{seq4}
\begin{equation*}\label{seq12}
\bigwedge^{c}\SN_{Y, \PP^N} =\bigwedge^{c-1} \mathcal{H}om_{Y}(\SI_{\WY,\PP^N}/\SI_{Y,\PP^N}^2, \SO_Y)\otimes \SE^{-1}.
\end{equation*}
So from~\eqref{seq10} we obtain~\eqref{SE'=SE-2}.
\end{Proof}
\begin{theorem}\label{hilbert}
Let $\WY$ be a projective $K3$ carpet on an Enriques surface $Y$ embedded in $\PP^{N}$ as in Theorem~\ref{embsmoothing}.
Then the Hilbert point of $\,\WY$ is nonsingular.
\end{theorem}

\begin{Proof}
We have proved in Theorem~\ref{embsmoothing} that $\WY$ admits an embedded smoothing.
Moreover, from Theorem~\ref{dim.19}, we know that for any $K3$ surface $X \subset \PP^N$ we have $H^1 (\SN_{X,\PP^N})= H^2 (\SN_{X,\PP^N})=0$.
So, from an straightforward computation, we see that the dimension of a component parametrizing $K3$ surfaces in $\PP^{N}$ is $18+2g(N+1)$. Therefore the $K3$ carpet $\WY$ represents a smooth point in the Hilbert scheme iff $h^0(\SN_{\WY,\PP^{N}})=18+2g(N+1)$.\\
As a consequence of Theorem~\ref{embsmoothing} or by direct computation using
the sequences in Lemma~\ref{Hilbert}, we see that the Euler characteristic is $\chi(\SN_{\WY,\PP^{N}})=18+2g(N+1)$.
Therefore we have to show that
\begin{equation*}
h^1(\SN_{\WY,\PP^{N}})-h^2(\SN_{\WY,\PP^{N}})=0.
\end{equation*}
Indeed, first we see at once that
\begin{equation*}\label{h-Y=0}
H^1(\SN_{Y,\PP^{N}})=H^2(\SN_{Y,\PP^{N}})=H^2(\SN_{Y,\PP^{N}}\otimes \omega_Y)=0.
\end{equation*}
In addition,~\eqref{dim6} says that
\begin{equation*}\label{h-Y=0.2}
H^1(\SN_{Y,\PP^{N}}\otimes \omega_Y)=0.
\end{equation*}
Therefore, from the sequences in Lemma~\ref{Hilbert}, we obtain
\begin{equation*}\label{for1}
H^2(\mathcal{H}om_{Y}(\SI_{\WY,\PP^{N}}/\SI_{Y,\PP^{N}}^2, \SO_Y))=0,
\end{equation*}
\begin{equation*}\label{for2}
H^2(\mathcal{H}om_{Y}(\SI_{\WY,\PP^{N}}/\SI_{Y,\PP^{N}}^2, \SO_Y)\otimes  \omega_Y)=0,
\end{equation*}
\begin{equation*}\label{for3}
H^1(\mathcal{H}om_{Y}(\SI_{\WY,\PP^{N}}/\SI_{Y,\PP^{N}}^2, \SO_Y))=H^2(\omega_Y^{-1})=\CC,
\end{equation*}
\begin{equation*}\label{for4}
H^1(\mathcal{H}om_{Y}(\SI_{\WY,\PP^{N}}/\SI_{Y,\PP^{N}}^2, \SO_Y)\otimes \omega_Y)=H^2(\SO_Y)=0.
\end{equation*}
Then we obtain
\begin{equation*}\label{for8}
H^1({\SN_{\WY,\PP^{N}}|}_Y \otimes \omega_Y)=0,
\end{equation*}
\begin{equation*}\label{for5}
H^2({\SN_{\WY,\PP^{N}}|}_Y \otimes \omega_Y)=H^2(\omega_Y^{-1})=\CC,
\end{equation*}
\begin{equation*}\label{for6}
H^2({\SN_{\WY,\PP^{N}}|}_Y)=0,
\end{equation*}
and
\begin{equation*}\label{for7}
H^1({\SN_{\WY,\PP^{N}}|}_Y)=0 \;\, \mathrm{or} \;\, \CC.
\end{equation*}
Finally, from sequence~\eqref{first-sequence-bis-twisted}, we see that $h^1(\SN_{\WY,\PP^{N}})-h^2(\SN_{\WY,\PP^{N}})=0 \;\,\mathrm{or}\;\, -1$, but now observe that, since our component has dimension $18+2g(N+1)$, we know that $h^0(\SN_{\WY,\PP^{N}})\geq \chi(\SN_{\WY,\PP^{N}})$, so $h^1(\SN_{\WY,\PP^{N}})-h^2(\SN_{\WY,\PP^{N}}) \geq 0$.
\end{Proof}

\providecommand{\bysame}{\leavevmode\hbox to3em{\hrulefill}\thinspace}
\providecommand{\MR}{\relax\ifhmode\unskip\space\fi MR }
\providecommand{\MRhref}[2]{%
  \href{http://www.ams.org/mathscinet-getitem?mr=#1}{#2}
}
\providecommand{\href}[2]{#2}


\begin{thebibliography}{BE95}

\bibitem[BE95]{BE95}
D. Bayer and D. Eisenbud, \emph{Ribbons and their canonical embeddings},
  Trans. Amer. Math. Soc. \textbf{347} (1995),  719--756.


\bibitem[Con00]{Conrad00}
B. Conrad, \emph{Grothendieck duality and base change}, Lecture Notes in
  Mathematics, vol. 1750, Springer-Verlag, Berlin, 2000.


\bibitem[Fon93]{Fong}
L.Y. Fong, \emph{Rational ribbons and deformation of hyperelliptic curves},
J. Algebraic Geom. {\bf 2} (1993),  295--307.


\bibitem[GGP05]{GGP}
F.J. Gallego, M. Gonz\'alez, and B.P. Purnaprajna,
\emph{Deformation of finite morphisms and smoothing of ropes}, preprint, arXiv:math.AG/0502467.

\bibitem[GP97]{GP}
F.J. Gallego and B.P. Purnaprajna, \emph{Degenerations of ${K}3$ surfaces in projective space},
Trans. Amer. Math. Soc. {\bf 349} (1997), 2477--2492.

\bibitem[Gon06]{Gon}
M. Gonz{\'a}lez, \emph{Smoothing of ribbons over curves},
J. reine angew. Math. {\bf 591} (2006), 201--235.

\bibitem[Gro61]{EGA3-1}
A. Grothendieck, \emph{{E}{G}{A} {I}{I}{I}, \'{E}tude cohomologique des
  faisceaux coh\'erents. (premi\`ere partie.)}, Publ. {M}ath. {I}{H}{E}{S},
  vol.~11, 1961.

\bibitem[Har66]{RD}
R. Hartshorne, \emph{Residues and duality}, Lecture notes of a seminar on
  the work of A. Grothendieck, given at Harvard 1963/64. With an appendix by P.
  Deligne. Lecture Notes in Mathematics, No. 20, Springer-Verlag, Berlin, 1966.

\bibitem[Har77]{Hartshorne77}
\bysame, \emph{Algebraic geometry}, Springer-Verlag, New York, 1977, Graduate
  Texts in Mathematics, No. 52.

\bibitem[Hor74]{Hor}
E. Horikawa, \emph{On deformations of holomorphic maps. II},
J. Math. Soc. Japan  \textbf{26}  (1974), 647--667.

\bibitem[HV85]{HV}
K. Hulek and A. Van de Ven, \emph{The Horrocks-Mumford bundle and the
  Ferrand construction},
  Manuscripta Math.  \textbf{50}  (1985), 313--335.


\bibitem[Kle80]{Kleiman80}
S.~L. Kleiman, \emph{Relative duality for quasicoherent sheaves},
  Compositio Math. \textbf{41} (1980), no.~1, 39--60.

\bibitem[Rei88]{Reider}
I. Reider, \emph{Vector bundles of rank $2$ and linear systems on
  algebraic surfaces},  Ann. of Math. (2)  \textbf{127}  (1988), 309--316.


\bibitem[SD74]{SaintDonat74}
B.~Saint-Donat, \emph{Projective models of {$K3$} surfaces}, Amer. J. Math.
  \textbf{96} (1974), 602--639.

\bibitem[Ser06]{Sernesi06}
E. Sernesi, \emph{Deformations of Algebraic Schemes}, Springer-Verlag, 2006,
Grundlehren der mathematischen Wissenschaften, Vol. 334.

\bibitem[Zar95]{Zar} O. Zariski, \emph{Algebraic surfaces}.
With appendices by S. S. Abhyankar, J. Lipman and D. Mumford.
Preface to the appendices by Mumford.
Reprint of the second (1971) edition.
Classics in Mathematics. Springer-Verlag, Berlin, 1995.


\end{thebibliography}
\end{document}